\documentclass[a4paper,11pt]{article}
\usepackage{amsmath}
\usepackage[T1]{fontenc}
\usepackage{amssymb}
\usepackage[english]{babel}
\usepackage{graphicx}

\ifx\JPicScale\undefined
\def\JPicScale{1}
\fi
\unitlength \JPicScale mm
\newtheorem{definition}{\textbf{\underline{Definition}}}
\newtheorem{lemma}{\textbf{\underline{Lemma}}}
\newtheorem{theorem}{\textbf{\underline{Theorem}}}
\newtheorem{proposition}{\textbf{\underline{Proposition}}}
\begin{document}

\begin{center}
{\Huge Work version of : existence of a translation invariant measure on $l^{\infty}$}

\vspace{10 mm}

{\LARGE Larrieu Jean-Yves}

{\Large E-mail: maths.larrieu@gmail.com}

\vspace{5 mm}
\end{center}

\begin{quotation}
{\Large Abstract } We describe two construction processes of
relevant measures, one in any non-empty compact metric space, and the other in the space $l^{\infty}(\mathbf{R})$. Both have invariance properties with respect to maps defined in a natural
way in these spaces. These properties imply that these measures are appropriate generalisations of the Lebesgue measure.
Results about their uniqueness are showed, and some applications and complementary properties are quickly studied.

{\Large R\'{e}sum\'{e} } Cet article d\'{e}crit deux proc\'{e}d\'{e}s de
construction de mesures pertinentes, l'un dans tout espace m%
\'{e}trique compact non vide, l'autre dans l'espace $l^{\infty}(\mathbf{R})$. Toutes deux poss\`{e}dent des propri\'{e}t\'{e}%
s d'invariance vis \`{a} vis d'applications naturellement d\'{e}finies dans
ces espaces. Ces propriétés montrent que ces mesures sont de bonnes généralisations de la mesure de Lebesgue. On prouve ensuite des théorèmes au sujet de leur unicité, puis on donne rapidement quelques applications de nos résultats.
\end{quotation}

\section*{Keywords:} Integral geometry, compact metric space, probability, measure invariant by isometries, infinite dimensional integration.

\section*{Introduction}

\section{Existence of an invariant by translations Borel measure on $l^{\infty} (\mathbf{R})$}

We shall now show our most impressive application of the geometry and probability theorem: on the space $l^{\infty}$ of bounded real sequences endowed with the topology given by the supremum norm, there is a measure defined on Borel sets which is invariant by translations. Moreover, this measure is locally finite and non zero.

\subsection{Construction of the measure on a unit cube}

\subsubsection{Construction of a weak measure invariant by translations in a unit cube}

Let $E = l^{\infty}$ the Banach space of bounded real sequences, endowed with the supremum norm $\Vert . \Vert_{\infty}$. Denote $(x_i)$ the canonical Schauder basis of $E$ (in weak topology), that is each $x_i$ is the sequence whose all values are zero, except the i-th which equals $1$; and denote $(f_i)$ the dual family of $(x_i)$: $f_i$ is the continuous form which gives the i-th coordinate of a vector. Let $\Vert .\Vert_w = \sum_i \frac{1}{2^i}\vert f_i(.)\vert$ be the weak norm. It defines a weaker topology than the usual one. We call the topology associated to $\Vert .\Vert_{\infty}$ the strong topology. Finally, denote $C = \bigcap_i f_i^{-1}([0;1])$ the unit cell of $E$ (it is the strong ball of $E$ centred in $\sum_i \frac{1}{2}x_i$ and with radius $1/2$).

$\Vert .\Vert_w$ is known to define the weak topology of $E$, that is the topology of the simple convergence of all coordinates. $C$ is compact for this topology. Then, our main theorem says that there is a measure $\mu_{C}$ on $C$ defined on weak Borel sets (that is, the set spanned as a $\sigma$-algebra by the weak topology of $C$). Moreover, this measure is invariant by partial $\Vert .\Vert_w$-isometries, peculiarly by translations defined on weak opens.

\vspace{2 mm} 

Now, let us now that the invariance by translations is in reality true for all weak Borel sets.

In order to prove that, we need the following definition.

\begin{definition}
A parallelepiped of $E$ is a set $\bigcap_i f_i^{-1}(I_i)$ where the $I_i$ are intervals of $\mathbf{R}$, or the empty set. If all the $I_i$ are closed, the parallelepiped is said closed, and if all the $I_i$ are open, the parallelepiped is said open.
\end{definition}

In fact, closed parallelepipeds are strongly and weakly closed, and the open ones are only strongly open.

\vspace{2 mm}

First of all, we show that the invariance is true for open parallelepipeds and for finite dimensional translations, that is translations whose vector $t$ is such as $t = \sum_{i\leq N} \lambda_i x_i$ for any $N$ and suitable $(\lambda_i)$ (set $\lambda_i = 0$ for $i > N$). Let then $B$ be an open parallelepiped of $E$ included in $C$ such as $B + t \subset C$. $B = \bigcap_i f_i^{-1}((c_i-r_i;c_i+r_i))$, where $(r_i)$ are positive real numbers and $(c_i)_i$ the center of $B$. One has $c_i-r_i > 0$ and $c_i+r_i < 1$ for all $i$, and since $B+t = \bigcap_i f_i^{-1}((c_i + \lambda_i-r_i;c_i+\lambda_i + r_i))$, $c_i+\lambda_i-r_i > 0$ and $c_i+\lambda_i+r_i < 1$ for all $i$. Choose a $\varepsilon >0$ such as $\varepsilon < \min(\lbrace c_i-r_i, 1-c_i-r_i, c_i+\lambda_i-r_i, 1-c_i-\lambda_i-r_i, i\leq N \rbrace)$. Then $C\cap\bigcap_{i\leq N} f_i^{-1}((c_i -r_i-\varepsilon;c_i + r_i + \varepsilon))$ and $C\cap\bigcap_{i\leq N} f_i^{-1}((c_i -r_i+\lambda_i-\varepsilon;c_i +\lambda_i+ r_i + \varepsilon))$ are weak opens of $C$ in relation by the translation of vector $t$. Hence the measures $\mu_{C}$ restricted to each of them are in relation by this translation, and $\mu_{C} (B) = \mu_{C} (B+t)$.

Note then that for any $f_i$, and for any $x \in (0,1)$, $f_i^{-1}(x)$ has measure $0$, since for all $y \in (0,1)$ and all $\alpha$ small enough, the $f_i^{-1}((y-\alpha,y+\alpha))$ are isometric for the weak distance. Hence $f_i^{-1}(x)$ in intrinsically null for $x \in (0,1)$, and even if $x \in [0,1]$. Since the boundary of any parallelepiped is a denumerable union of sets included in such $f_i^{-1}(x)$, it is null, and the measures of closed parallelepipeds are invariant by finite dimensional translations too.

Now,we show that the invariance is true for parallelepipeds and for any translation.

Let $B$ be a closed parallelepiped and $t = \sum_i \lambda_i x_i$ such that $B$ and $B+t$ are contained in $C$. Denote $B_N = B + t_N$ where $t_N = \sum_{i\leq N} \lambda_i x_i$. $(t_N)$ weakly converges to $t$. Note that $B_N \subset C$ for all $n$. The space of weakly closed subsets (as closed parallelepipeds for example) of $C$, endowed with the Hausdorff distance associated with the weak distance is well known to be a compact set. Thus, the following proposition holds :

\vspace{2 mm}

\begin{proposition}
Let $K$ be a compact metric space, $[K]$ be the set of closed of $K$ endowed with the Hausdorff distance, and $m$ a probability defined on Borel sets of $K$. If $(A_n)$ is a converging sequence of $[K]$, one has: 
$$ \limsup_{n \rightarrow +\infty} m(A_n) \leq m(\lim_{n \rightarrow +\infty} A_n)$$
\end{proposition}

\vspace{2 mm}

Let us show that $B_N$ approaches $B$ in Hausdorff weak distance. Choose a $\varepsilon > 0$, and $N$ such that $\sum_{i>N} \frac{1}{2^i} < \varepsilon$. Then, for any $p \geq N$, the Hausdorff weak distance between $B + t$ and $B_p$ is smaller than $\varepsilon$. Indeed, for any $x$ in $B + t$, $x - \sum_{i> p} \lambda_i x_i$ is in $B_p$, and $\Vert \sum_{i> p} \lambda_i x_i \Vert_w < \varepsilon$. Hence $\inf_{y \in B_p} \Vert y-x\Vert_w <\varepsilon$. We can show the same way that $\inf_{x \in B} \Vert y-x\Vert_w <\varepsilon$ for any $y$ in $B$. Hence $\limsup_{n \rightarrow +\infty} \mu_{C} (B_N) = \mu_{C} (B) \leq \mu_{C}(B+t)$ for any closed parallelepiped $B$, and any vector $t$ such that the partial translation is well defined.

Hence, $\mu_{C} (B+t)\leq \mu_{C} (B+t-t) = \mu_{C}(B)$, and parallelepipeds keep their measure by translation.

\vspace{2 mm}

Now, it is clear that the family of parallelepipeds is closed under finite intersections, and that it contains a basis of the weak topology (indeed, it contains the $f_i^{-1}((a,b))\cap C$ for all $i$ and all $(a,b) \subset \mathbf{R}$ whose set spans the weak topology). Consider then $t$ a vector of $C$, and the family of weak Borel sets $S$ of $C$ such that $S + t$ is in $C$. The biggest one is $B=C \cap (C-t)$, and it is a parallelepiped. Any parallelepiped of $B$ is measure-invariant by the $t$-translation, and the family of Borel sets of $B$ measure-invariant by this translation is a Dynkin system. Hence it contains the $\sigma$-algebra generated by the family of parallelepipeds of $B$, and consequently the $\sigma$-algebra generated by the weak topology of $B$ (which is weakly separable). So, all weak Borel sets are measure-invariant by partial translations.

\subsubsection{Construction of the measure for all Borel sets}

We will now use a method which looks like our construction of the metrically compatible probability. The lemma 13 is the translation of the proposition 4, the lemma 14 corresponds to lemmas 3 and 5, and lemma 15 to lemma 7. The present construction is not exactly the same as before since the situation is a bit different, but the underlying ideas are the same. We begin with definitions and proofs imitating what happens with the algebra of sets $M$.

\vspace{2 mm}

\begin{definition}
Let $E$ be a set. A family $F$ of subsets of $E$ stable by countable union and containing the empty set is called a $u$-system.

Given a family $G$ of subsets of $E$, the smallest $u$-system containing $G$ is called the $u$-system generated by $G$.
\end{definition}

\vspace{2 mm}

Note that the $u$-system generated by any family of sets always exist: just take the intersection of the $u$-systems containing the family.

With the notations of the paragraph, if $A$ is a subset of $C$, we denote $u (A,C)$ the $u$-system generated by the closed parallelepipeds of $C$ included in $A$. Obviously, for any $A$, $u (A,C)$ is contained in the weak Borel $\sigma$-algebra of $C$. Moreover, we have:

\vspace{2 mm}

\begin{lemma}
Let $O$ and $O'$ be disjoint strong opens of $C$. Then: 
$$u (O\amalg O',C) = \lbrace b\amalg b', b\in u (O,C), b'\in u (O',C)\rbrace.$$ Consequently, for any $B \in u (O\amalg O',C)$, $B\cap O \in u (O,C)$ and $B\cap O'\in u (O',C)$.
\end{lemma}

\underline{\textbf{Demonstration:}} On one hand, it is obvious that $u (O,C) \subset u (O\amalg O',C)$ and that $u (O',C) \subset u (O\amalg O',C)$. Hence $\lbrace b\amalg b', b\in u (O,C), b'\in u (O',C)\rbrace \subset u (O\amalg O',C)$ thanks to stability by unions.

On the other hand, $\lbrace b\amalg b', b\in u (O,C), b'\in u (O',C)\rbrace$ contains all parallelepipeds included in $O\amalg O'$. Indeed, since parallelepipeds are connected, if a parallelepiped is included in $O\amalg O'$, it is included in either $O$ or $O'$.

We just have to prove that $\lbrace b\amalg b', b\in u (O,C), b'\in u (O',C)\rbrace$ is a $u$-system. But it is clear that it contains $\emptyset$, and that it is stable by countable unions. 

Consequently, $\lbrace b\amalg b', b\in u (O,C), b'\in u (O',C)\rbrace$  is a $u$-system containing parallelepipeds included in $O\amalg O'$, hence: $u (O\amalg O',C) \subset \lbrace b\amalg b', b\in u (O,C), b'\in u (O',C)\rbrace$, which leads to the conclusion.

For the last point, if $B\in u (O\amalg O',C)$, $B$ can be written $B = b \cup b'$ with $b \in u (O,C)$ and $b' \in u (O',C)$, and $B\cap O = b$, $B\cap O' = b'$.

\begin{flushright}
$\blacksquare$
\end{flushright}

Now, we define a new set function in order to define the desired measure. For all subsets $A$ of $C$, 

$$ m'(A) = \sup_{B \in u (A,C)} \mu_{C} (B).$$

$m'$ is an increasing function, $m'(\emptyset ) = 0$,  and $m'$ is invariant by well defined partial translations. Moreover, if $O$ is a weak open, $m'(O) = \mu (O)$. Indeed, $O$ can be written as a denumerable union of weak closed positively included in $O$, and each of them can be covered by a finite union of weak balls positively included (for the weak norm) in $O$. Just use the weak compactness of $C$ to show this, taking a finite sub-covering of a covering by small enough weakly open balls. Hence $O$ can be written as a denumerable union of weakly closed weak balls, and these are closed parallelepipeds. Hence $O \in u (O,C)$, and $m'(O) = \mu_{C} (O)$.

\vspace{2 mm}

Next, if $(A_i)$ is a sequence of disjoint subsets of $C$, for any $\varepsilon > 0$, choose a $B_i$ in $u (A_i,C)$ such that $m'(A_i) \geq \mu_{C} (B_i) \geq m'(A_i) - \frac{\varepsilon}{2^{i+1}}$. Then $m'(\coprod_i A_i) \geq \mu_{C} (\coprod_i B_i) \geq \sum_i m'(A_i) - \varepsilon$. Consequently, $m'(\coprod_i A_i) \geq \sum_i m'(A_i)$.

\vspace{2 mm}

Now, if $O_i$ are disjoint strong opens of $C$, we have: 
$$m'(\coprod_i O_i) \geq \sum_i m'(O_i).$$
Choose then any $\varepsilon >0$, and a $B \in u(\coprod_i O_i,C)$ such that $m'(\coprod_i O_i) \leq \mu_{C} (B) + \varepsilon$. With the lemma, $B\cap O_i \in u(O_i,C)$ for all $i$. Indeed, $O_i$ and $\coprod_{j\neq i} O_j$ are disjoint opens. Hence: $m'(\coprod_i O_i) \leq \mu_{C} (B) + \varepsilon = \sum_i \mu_{C} (B \cap O_i)\leq \sum_i m'(O_i)$. So: $m'(\coprod_i O_i) = \sum_i m'(O_i)$. $m'$ appears to behave like a measure for open sets. This fact is quite encouraging for our purpose. Now, let us show :

\vspace{2 mm}

\begin{lemma}
$m'(\overline{O})= m'(O)$ for any strong open $O$ of $C$.
\end{lemma}

\underline{\textbf{Demonstration:}} If $O$ is a strong open of $C$, we have $m'(\overline{O}) \geq m'(O)$, where $\overline{O}$ denotes the strong closure of $O$. And for any $\varepsilon > 0$, $m'(\overline{O}) \leq \mu_{C} (\bigcup_i B_i) + \varepsilon$ for some denumerable union $\bigcup_i B_i$ of closed parallelepipeds $B_i$ of $\overline{O}$. Now, since $\mu_{C}$ is regular, for any $i$, there is a closed parallelepiped $C_i$ positively included in $B_i$ (for the strong norm) such that $\mu_{C} (B_i) \leq \mu_{C} (C_i) + \frac{\varepsilon}{2^{i+1}}$. Hence $\mu_{C} (\bigcup_i B_i \setminus \bigcup_i C_i) \leq \sum_i \frac{\varepsilon}{2^{i+1}} = \varepsilon$, and $C_i \subset O$ for all $i$. So: 
$m'(\overline{O}) \leq \mu_{C} (\bigcup_i B_i) + \varepsilon \leq \mu_{C} (\bigcup_i C_i) + 2\varepsilon \leq m'(O) + 2\varepsilon $. This is true for all $\varepsilon >0$, so: $m'(\overline{O}) = m'(O)$ for all strong opens.
\begin{flushright}
$\blacksquare$
\end{flushright}

\vspace{2 mm}

Following the proof of our first main theorem, we show now :

\vspace{2 mm}

\begin{lemma}
Let $(O_i)$ be a countable family of strong opens, not assumed to be disjoint. Then :

$$\sum_i m'(O_i)\geq m'(\bigcup_i O_i).$$
\end{lemma}

\underline{\textbf{Demonstration:}}
$$ \sum_i m'(O_i)\geq \sum_i m'(O_i\setminus \overline{\bigcup_{j<i} O_j})=m'(\coprod_i O_i\setminus \overline{\bigcup_{j<i} O_j})$$
$$= m'(\overline{\coprod_i O_i\setminus \overline{\bigcup_{j<i} O_j}}) \geq m'(\bigcup_i (\overline{ O_i }\setminus\overbrace{\overline{\bigcup_{j<i} O_j}}^{\circ} ))$$
$$\geq m'(\bigcup_i (\overline{ O_i }\setminus \bigcup_{j<i} \overline{O_j} ))\geq m'(\bigcup_i \overline{ O_i })\geq m'(\bigcup_i O_i ).$$
\begin{flushright}
$\blacksquare$
\end{flushright}

\vspace{2 mm}

The construction is nearly finished. Let us define: 

$$m(A) = \inf_{A\subset O \mbox{ strong open}} m'(O).$$

$m$ is increasing, $m(\emptyset ) = 0$, and $m$ equals $m'$ on strong opens. Consequently, $m$ equals $\mu_{C}$ on weak opens, and $m$ is additive on disjoint strong opens. 

Let us now show that $m$ is an exterior measure: let $(A_i)$ be a countable family of sets. Choose a $\varepsilon > 0$. For all $i$, $m(A_i) \geq m'(O_i) - \frac{\varepsilon}{2^{i+1}}$ for some strong open $O_i \supset A_i$. Then $m(\bigcup_i A_i ) \leq m' (\bigcup_i O_i) \leq \sum_i m'(O_i)\leq \sum_i m(A_i)+\varepsilon$. This is true for all $\varepsilon > 0$, hence $m$ is an exterior measure.

Therefore, with the theorem of metric exterior measures, $m$ is a measure defined on strong Borel sets of $C$, and $m$ prolongs $\mu_{C}$. Consequently, it is a probability and it is outer regular and inner regular by closed sets, since $C$ is a metric space (see \cite{Laudenbach} for this theorem). We use the following definition:

\vspace{2 mm}

\begin{definition} 
A measure is said inner regular by closed sets if the measure of any measurable set can be approached as close as wished by the measure of a closed contained in.
\end{definition}

\vspace{2 mm}

A question naturally arises now : is $m$ the restriction to strong Borel sets of the completion of $\mu_{C}$ by null sets ? A positive answer would really help for the comprehension of $m$. In fact, it is true. To show this, let $A$ be a strong closed. For all positive integer $n$, there is a strong open $O_n\supset A$ such as $m'(O_n)\geq m(A)\geq m'(O_n) - \frac{1}{2^n}$, and we can choose the $O_n$ in order to have $A= \bigcap_n O_n$. For all $k \in \mathbf{N}^*$, $A \subset \bigcap_{n\leq k} O_n$ and $m'(\bigcap_{n\leq k} O_n)\geq m(A)\geq m'(\bigcap_{n\leq k} O_n) - \frac{1}{2^k}.$ So $\lim_{k \rightarrow \infty}m'(\bigcap_{n\leq k} O_n)= m(A)$.

Now, if $U$ is any set, one can choose a $V \in u(U,C)$ such that : $\mu_C \leq m'(U) \leq\mu_C(V) + \varepsilon/2$ for any $\varepsilon >0$. And since $V$ is a countable union of closed parallelepipeds, we can take a finite union $W$ of parallelepipeds such that $W \subset V$ and $\mu_C(V)-\mu_C(W)<\varepsilon/2$. Then $\mu_C(W)\leq m'(V) \leq \mu_C(W)+\varepsilon$ and $W$ is a weak closed.

Then, choose a $\varepsilon >0$. With a recursive process, for $k = 1$, take a weak closed $F_1 \subset \bigcap_{n\leq 1} O_n$ such as $\mu_C (F_1)\leq m'(\bigcap_{n\leq 1} O_n)\leq \mu_C (F_1) + \varepsilon$.

Next, for $k=2$, we have $\bigcap_{n\leq 2} O_n \setminus F_1 \subset \bigcap_{n\leq 1} O_n \setminus F_1$, hence it is possible to choose a weak closed $F_2 \subset F_1\cap \bigcap_{n\leq 2} O_n$ such that $\mu_C (F_2)\leq m'(F_1\cap\bigcap_{n\leq 2} O_n)\leq \mu_C (F_2) + \varepsilon/2$, so with $\mu_C(F_2) \leq m'(\bigcap_{n\leq 2} O_n)\leq m'(\bigcap_{n\leq 2} O_n\setminus F_1)+m'(\bigcap_{n\leq 1} O_n)$. Indeed, both $\bigcap_{n\leq 1} O_n$ and $\bigcap_{n\leq 2} O_n\setminus F_1$ are strong opens, and $\bigcap_{n\leq 2} O_n \subset \bigcap_{n\leq 1} O_n \cup (\bigcap_{n\leq 2}O_n\setminus F_1).$ So $\mu_C(F_2) \leq m'(\bigcap_{n\leq 2} O_n) \leq \mu_C(F_2) + \varepsilon + \varepsilon /2.$

We build this way a decreasing sequence $F_k$ of weak closed sets such that for all $k$, $F_k \subset \bigcap_{n\leq k} O_n$ and $\mu_C (F_k)\leq m'(\bigcap_{n\leq k} O_n)\leq \mu_C (F_k) + \varepsilon \sum_{i=0}^{k-1}\frac{1}{2^i}$.

Therefore $\bigcap_n F_n \subset \bigcap_n O_n$ and for all $k$, $$\mu_C (\bigcap_n F_n)\leq m'(\bigcap_n O_n)\leq \lim_{k \rightarrow \infty}m'(\bigcap_{n\leq k} O_n)\leq \mu_C (F_k) + 2\varepsilon ,$$ hence $\mu_C (\bigcap_n F_n)\leq m(A)\leq \mu_C (\bigcap_n F_n) + 2\varepsilon$, and $\bigcap_n F_n \subset \bigcap_n O_n = A$.

So, for all $\varepsilon > 0$, there is a weak closed $W \subset A$ such that $m(W)\leq m(A)\leq m(W) + \varepsilon$. Hence, any strong closed $A$ can be written as a union of a weak Borel set and a null set.

Hence the $\sigma$-algebra generated by the weak Borel sets and the $m$-null sets contains strong closed sets, so it contains the $\sigma$-algebra generated by strong closed, and thus strong Borel sets.

Now, let us show that $m$-null sets are $\mu_C$-null sets. Let $N$ be any $m$-null set. For all $\varepsilon > 0$, there is a strong open $O$ such that $N \subset O$ and $m(O) < \varepsilon$, since $m$ is regular. But $m(O) = m'(O)$, and $m'(O\amalg \complement_C \overline{O}) = m'(O)+ m'(\complement_C \overline{O})=m'(\overline{O\amalg \complement_C \overline{O}})=m'(C)=1.$ Thus, $m'(O) = 1 - m'(\complement_C \overline{O})< \varepsilon,$ and $m'(\complement_C \overline{O}) > 1 - \varepsilon.$ So, with the definition of $m'$, there is a countable union of closed parallelepipeds, hence a weak Borel set $B$, such that $B \subset \complement_C \overline{O}$, and $\mu_C(B) > 1- 2\varepsilon$. Then, $N \subset O \subset \overline{O} \subset \complement_C B$, and $\mu_C (\complement_C B) < 2\varepsilon.$ Since for all $\varepsilon > 0$, $\complement_C B$ is a weak Borel set, $N$ is $\mu_C$-null. 

So:
\vspace{2 mm}

\begin{theorem} 
The completion of $\mu_C$ on $C$ by null sets is defined on strong Borel sets.
\end{theorem}

\vspace{2 mm}

Last point, it is obvious that $m$, the restriction of the completion of $\mu_C$ to strong Borel sets, is invariant by translations, since both $\mu_C$ and the family of null sets are invariant by translations.

\subsubsection{A patching measures theorem}

The further step of our construction uses the following theorem.

\begin{theorem} 
Let $E$ be a set endowed with a $\sigma$-algebra $M$ such that there is a family $(B_i)_{i \in I}$ of measurable sets covering $E$, and such that for all $i\in I$ there is a non-negative measure $m_i$ on $B_i$ defined on the induced $\sigma$-algebra. There is then a measure $\mu$ on $E$ defined on $M$, canonically chosen, such as: 
$$\mu (B) = \sup \sum_{i} m_i(A_i),$$
where the supremum is taken over disjoint families $(A_i)_{i \in I}$ of measurable sets contained in $B$, with for all $i$, $A_i \subset B_i$, and such that the set of $i$ satisfying $A_i \neq \emptyset$ is countable. It is called the supremum measure of the $(m_i)_{i \in I}$.

Moreover, assuming that for all $i$ and $j$ in $I$, and all measurable set $B$ of $M$ contained in $B_i\cap B_j$, $m_i(B) = m_j(B)$ then for all $i\in I$, and all measurable set $B$ contained in $B_i$, $\mu (B) = m_i(B)$. This last measure is called the canonical patching measure of the family $(m_i)_{i \in I}$.
\end{theorem}

\vspace{2 mm}

The fact that no cardinal assumption has been made about $I$ has to be pointed out.

\vspace{2 mm}

\underline{\textbf{Demonstration:}} Let us use the definition of the theorem: for all measurable set $B$ of $E$, 
$$\mu (B) = \sup \sum_{i} m_i(A_i),$$
where the supremum is taken over disjoint families $(A_i)_{i \in I}$ of measurable sets contained in $B$, with for all $i$, $A_i \subset B_i$, and such that the set of $i$ satisfying $A_i \neq \emptyset$ is countable.

Let us show that $\mu$ is a measure. Let $(C_j)_{j \in \mathbf{N}}$ be a countable family of disjoint measurable sets, and $C$ its union. We have to prove that $\mu(C) = \sum_{j} \mu (C_j)$.

Choose a $\varepsilon >0$. There is a family $(A_i)_{i \in I}$ of measurable sets contained in $C$, with for all $i$, $A_i \subset B_i$, such that the set of $i$ satisfying $A_i \neq \emptyset$ is countable, and such as: $\mu(C) \geq \sum_{i} m_i(A_i) \geq \mu(C) - \varepsilon.$ Then: $\mu(C) - \varepsilon \leq \sum_{i} \sum_{j} m_i(A_i \cap C_j) =  \sum_{j} \sum_{i} m_i(A_i \cap C_j) \leq \sum_{j} \mu(C_j)$. Indeed, one can swap the sum symbols in series whose terms are non-negative, and the $(A_i \cap C_j)_{i \in I}$ are disjoint measurable sets contained in $C_j$, with for all $i$, $A_i \cap C_j \subset B_i$, and such that the set of the $i$ satisfying $A_i\cap C_j \neq \emptyset$ is countable. This being true for all $\varepsilon >0$, we get an inequality.

Then, a $\varepsilon >0$ still being chosen, for all $j$, there is a disjoint family $(C_{i,j})_{i \in I}$ of measurable contained in $C_j$, with for all $i$, $C_{i,j} \subset B_i$, such that the set of $i$ satisfying $C_{i,j} \neq \emptyset$ is countable, and such as: $\mu (C_j) \leq \sum_i m_i(C_{i,j}) + \frac{\varepsilon}{2^j}.$ Then: $ \sum_j \mu (C_j) \leq 2\varepsilon + \sum_j \sum_i m_i(C_{i,j}) = 2\varepsilon + \sum_i \sum_j m_i(C_{i,j}) = 2\varepsilon + \sum_i m_i(\coprod_j C_{i,j}).$ Indeed, the $(C_{i,j})_{j \in \mathbf{N}}$ are disjoint and contained in $B_i$ for all $i$. Then: $ \sum_j \mu (C_j) \leq 2\varepsilon + \sum_i m_i(\coprod_j C_{i,j}) \leq \mu (C) + 2\varepsilon,$ since the $(\coprod_j C_{i,j})_{i \in I}$ are disjoint measurable sets contained in $C$, such that for all $i$, $\coprod_j C_{i,j} \subset B_i$, and such that the set of $i$ satisfying $\coprod_j C_{i,j} \neq \emptyset$ is countable. Indeed, $\coprod_j C_{i,j} \neq \emptyset \Leftrightarrow \exists j \in \mathbf{N} / C_{i,j} \neq \emptyset$, hence the set of $i$ such that $\coprod_j C_{i,j} \neq \emptyset$ is a countable union of countable sets. This being true for all $\varepsilon >0$, we get the reverse inequality.

This shows that $\mu$ is a measure. Let us prove now that it prolongs the $(m_j)$ under the last assumption: choose a $j$ in $I$, and $B$ a measurable set of $B_j$. First of all, $\mu(B) \geq m_j(B)$ (look at the definition). Then, for all $\varepsilon >0$, there is a disjoint family $A_i$ of measurable sets contained in $B$, such that for all $i$, $A_i \subset B_i$, and such that the set of $i$ satisfying $A_i \neq \emptyset$ is countable, with:  $\mu(B) \leq \sum_{i} m_i(A_i) + \varepsilon.$ But all the $A_i$ are contained in $B_j$, and with the theorem assumption: $m_i(A_i) = m_j(A_i)$. Hence, with $\sigma$-additivity, $$\mu(B) \leq m_j(\coprod_{i} A_i) + \varepsilon \leq m_j(B) + \varepsilon.$$ This being true for all $\varepsilon >0$, we get the conclusion.

Finally, this construction is canonical, since it uses no arbitrary choice.\begin{flushright}
$\blacksquare$
\end{flushright}

Note that when $I$ is countable, that patching measure is unique: to prove it, just build a countable measurable partition of $E$ whose elements are contained in the $B_i$. By instance, choose $P_n = B_n\setminus (\bigcup_{k\leq n-1}B_k).$ One can then show the uniqueness proving that for a measurable set $B$, the only possible measure $\mu(B)$ is: $\mu(B) = \sum_n m_n(B \cap P_n).$

At last, note that when $I$ is not countable, many measures can be compatible with restrictions. Indeed, chose $E = \mathbf{R}$ and $B_i = \lbrace i\rbrace$ for all $i \in \mathbf{R}$, and $m_i = 0$. Then the Lebesgue measure and the zero measure are compatibles with this system.

\subsubsection{Construction of the measure on $l^{\infty}$}

We use the previous theorem and the measure $m$ built before to construct our measure on the whole space $E = l^{\infty}$. Let $C_t = C +t$ for all $t \in E$ and $m_t$ be the translation of $m$ of vector $-t$. It is a probability defined on Borel sets of $C_t$, invariant by well defined partial translations. Check now that the family $m_t$ can be patched. Let $B$ be any Borel set of $C_t \cap C_{t'}$. We have: $B-t \subset C$ and $B-t' \subset C$. Since $B-t = B-t'+(t'-t)$, $m(B-t) = m(B-t')$, and $m_t (B) = m_{t'} (B)$, thanks to the invariance of $m$ by partial translations. The hypothesis of the patching measure theorem are satisfied, so the $(m_t)$ define a patching measure that we will denote $\mu$ on the $\sigma$-algebra of (strong) Borel sets of $E$. For any $t \in E$, it equals $m_t$, so for any strong ball $D$ of radius $1/2$, $\mu (D) = 1$, and $\mu $ is strongly locally finite. 

Moreover, looking at the definition of the patching measure, we see that $\mu $ is inner regular by weakly closed sets. Indeed, let $A$ be a Borel set of $E$. If $\mu (A)$ is finite, for any $\varepsilon >0$, one can choose $(A_{t_n})$, a sequence of disjoint Borel sets included for all $n$ in $(C+t_n )\cap A$ such that $\mu (A) \leq \sum_{n} m_{t_n}(A_{t_n}) + \varepsilon$. We can then take an integer $N$ such that $\mu (A) \leq \sum_{n\leq N} m_{t_n}(A_{t_n}) + 2\varepsilon$, and for each $n \leq N$, a weak closed $F_n$ included in $A_{t_n}$ such that $m_{t_n}(A_{t_n}) \leq m_{t_n}(F_n) + \frac{\varepsilon}{2^{n+1}}$. Then: $\mu (A) \leq \sum_{n\leq N} m_{t_n}(F_n) + 3\varepsilon = \mu (\coprod_{n\leq N} F_n) + 3\varepsilon $, and $\coprod_{n\leq N} F_n$ is a weak closed contained in $A$. If $\mu (A)$ is infinite, the same kind of reasoning gives the result.

Finally, $\mu $ is invariant by translations. Indeed, let $A$ be any Borel set, and $t$ be any vector of $E$. Assume $\mu (A+t)$ is finite. Then for all $\varepsilon > 0$, $\mu (A+t) \geq \sum_{i} m_{t_i}(A_i)\geq \mu (A+t) - \varepsilon $ for a denumerable family of disjoint Borel sets $A_i$ included in $(C+t_i)\cap (A+t)$. And: $\sum_{i} m_{t_i}(A_i) = \sum_{i} m(A_i - t_i) = \sum_{i} m_{-t+t_i}(A_i -t)$. The $(A_i - t)$ are disjoint Borel sets included in $(C+t_i -t)\cap A$ for all $i$. Hence: $\mu (A+t) \leq \sum_{i} m_{t_i}(A_i) + \varepsilon \leq \sum_{i} m_{-t+t_i}(A_i -t) + \varepsilon \leq \mu (A)+ \varepsilon $ for all $A$ and $t$. Therefore, $\mu (A+t) \leq \mu (A)$. If $\mu (A+t)$ is infinite, we show similarly that: $\mu (A+t) \leq \mu (A)$. 

Then $\mu (A) \leq \mu (A+t-t) \leq \mu (A+t)$, and $\mu (A) = \mu (A+t)$ for all Borel set $A$ and all $t\in E$.

We showed: 

\vspace{2 mm}

\begin{theorem} 
There is a measure $\mu $ on $l^{\infty}$ defined on Borel sets which is locally finite, invariant by translations, and inner regular by weakly closed sets. More precisely $\mu (D)=1$ for any ball $D$ of radius $1/2$, and the restriction of $\mu $ to Borel sets of $D$ is weakly outer regular. 
\end{theorem}

\vspace{2 mm}

The construction could have been done for any non-empty parallelepiped instead of $C$.

Note that our theorem is not contradictory with the known theorems affirming that many obstructions prevent the existence of relevant measures in infinite dimensional spaces. Recall for example the following theorem (see \cite{sudakov}).

\vspace{2 mm}

\begin{theorem} 
On a locally convex, infinite dimensional, topological vector space, there does not exist any Borel 
measure which is $\sigma$-finite and invariant by all translations.
\end{theorem}

\vspace{2 mm}

In fact our measure is not $\sigma$-finite. This is due to the non-separability of $l^{\infty}$.

\vspace{1 mm}

In addition of this theorem, we have a characterisation of null sets for $\mu$ and a description of strong Borel sets of finite measure.

\vspace{2 mm}

\begin{theorem} 
The null sets for the measure $\mu$ previously defined on $l^{\infty}$ are exactly the sets $A$ such that for any ball $C$ of radius $1/2$, $\mu (C\cap A) =0$ (still denoting $\mu$ for its own completion by null sets).
\end{theorem}

\vspace{2 mm}

\underline{\textbf{Demonstration:}} It is obvious that null sets satisfy the condition. Conversely, if a set $A$ satisfies the condition, just look at the definition of $\mu$ as a patching measure. \begin{flushright}
$\blacksquare$
\end{flushright}

\begin{theorem} 
The strong Borel sets of $E$ of finite measure are exactly the unions of null sets with increasing countable unions of weakly closed sets.
\end{theorem}

\vspace{2 mm}

\underline{\textbf{Demonstration:}} Just use the fact that $\mu$ is inner regular by weakly closed sets. \begin{flushright}
$\blacksquare$
\end{flushright}

Finally, note the following fact. Let $C$ be the unit cell of $E$ associated with the canonical basis. One obviously has $\mu (C) =1$. But, denoting $L$ the lattice associated with the canonical basis $(x_i)$ (that is $L = \lbrace \sum_i \lambda_i x_i, (\lambda_i) \mbox{ bounded sequence of integers}\rbrace$), and $L' = \sum_i \frac{1}{2} x_i + L$ a translated lattice, we have the following property: for any $z \in L'$, $\mu ((z+C)\cap C)=0 $. Indeed, $\mu (z+C)\cap C \leq \mu (\bigcap_{n} f_n^{-1}([0;1/2])$ with invariance by translations ($(f_n)$ is the dual basis of $(x_i)$), and $\mu (\bigcap_{n\leq N} f_n^{-1}([0;1/2]) \leq (\frac{1}{2})^{N}$ for all $N \in \mathbf{N}$. Hence, a set can be non-null but can appear to be null along all the cells of a lattice. Consequently, the previous characterisation of null sets seems sharp. Peculiarly, knowing the restrictions of $\mu $ along all cells of the lattice associated with the canonical basis is not sufficient for knowing $\mu$ on the whole space $E$.

\subsection{Uniqueness of the measure}

Now, we give sufficient conditions for uniqueness of our measure.

\vspace{2 mm}

\begin{definition} 
Let $E$ be a metric measured space. A measurable set $A$ such that for any ball $B$ of radius $r>0$, $A\cap B$ is null is called a null set for $r$-balls.
\end{definition}

\vspace{2 mm}

\begin{theorem} 
In $l^{\infty}$, denote $C$ a ball of radius $1/2$. If a Borel measure $m$ is invariant by translations, with $m(C)=1$, and if any null set for $1/2$-balls is null, then $m$ equals the previously defined measure $\mu$.
\end{theorem}

\vspace{2 mm}

\underline{\textbf{Demonstration:}} Let $n$ be any natural number and $(x_i)$ the canonical basis of $E = l^{\infty}$. Consider the restriction of $m$ to the $\sigma$-algebra $B_n$ of Borel sets of $E$ which can be written $A\times Vect(x_i, i>n)$, where $Vect(x_i, i>n) = \lbrace \sum_{i>n} \lambda_i x_i, \lambda_i \in \mathbf{R}, (\lambda_i) \mbox{ bounded}\rbrace$. The push-forward measure of that restriction by $\Phi_n : E \rightarrow \mathbf{R}^n$ defined by: $\Phi_n (\sum_i \lambda_i x_i) = \sum_{i\leq n} \lambda_i e_i$ ($(e_i)$ being the canonical basis of $\mathbf{R}^n$) is just the Lebesgue measure on $\mathbf{R}^n$. Indeed, this measure is invariant by translations and the measure of the unit cell is $1$. These conditions are satisfied by the push-forward measure of $\mu $ by $\Phi_n$ too. Hence, since they characterize the Lebesgue measure, $\Phi_{n*}(m) = \Phi_{n*}(\mu)$, and $\mu $ and $m$ equal on $\bigcup_n B_n$. Indeed, for each $n$, $B_n$ is the push-forward $\sigma$-algebra of the Borel $\sigma$-algebra of $\mathbf{R}^n$ by $\Phi_n$. Therefore, they equal on weak opens of $C$, hence on the Dynkin system generated by weak opens, that is on the weak Borel $\sigma$-algebra of $C$. Consequently, their completions equal, and $\mu = m$ on the strong Borel sets of $C$. With invariance by translations, this true for any ball $B$ of radius $1/2$.

It follows that null sets of $E$ for $\mu$ and for $m$ are the same, using the third hypothesis of the theorem.

\vspace{2 mm}

Now, let us study the case of general strong Borel sets. Let $A$ be such a set.

With the notation of the following lemma, assume first that for a $\delta >0$ and a $r \in C$, $NZ_{\delta, > }(r)(A)= \lbrace z \in L / m((A-r)\cap (C+z))> \delta \rbrace$ is infinite. Note that with invariance by translations, $NZ_{\delta, > }(r)(A)= \lbrace z \in L / m(A\cap (C+z+r))> \delta \rbrace$. $A$ contains then a Borel set with measure as tall as wished for $m$ and $\mu$ (take unions of $A\cap \overbrace{C+z+r}^{\circ}$ which are disjoint). Then $m (A) = \mu (A) = +\infty$.

Assume now that $NZ_{\delta, > }(r)(A)$ is always finite. Since $m$ is a Borel measure on $E$ which is inner and outer regular for weak topology on every $1/2$-ball, with the following lemma, $A$ can be written $A = N\cup M$, with $N$ null for $1/2$-balls, hence null, and $M$ which can be covered by a countable family of $1/2$-balls, say $(B_n)$. Then : $m (A) = m(M) = \sum_n m(M\cap (B_n\setminus \bigcup_{k<n} B_k)) = \sum_n \mu (M\cap (B_n\setminus \bigcup_{k<n} B_k)) = \mu (M)= \mu (A)$. So $m= \mu $.
\begin{flushright}
$\blacksquare$
\end{flushright}

\vspace{2 mm}

\begin{definition} 
A set $A\subset E$ is said $\sigma$-finite for $\mu$ if there is a countable family of $1/2$-balls $(B_i)$ such that $A \setminus \bigcup_i B_i$ is $\mu$-null. Such a family $(B_i)$ is called an almost everywhere countable covering by $1/2$-balls.
\end{definition}

\vspace{2 mm}

\begin{lemma}[Key Lemma for $\sigma$-finiteness]
Let $(x_i)$ be the canonical basis of $E = l^{\infty}$, $L$ the associated lattice:

$$L = \lbrace \sum_i \lambda_i x_i, (\lambda_i) \mbox{ bounded sequence of integers}\rbrace$$

and $C = \lbrace \sum_i \lambda_i x_i, (\lambda_i) \in [0,1]^{\mathbf{N}}\rbrace$ its unit cell. Let $A$ be any strong Borel set of $E$, and denote $NZ_{\delta , > }(r)(A) = \lbrace z \in L / \mu ((A-r)\cap (C+z))> \delta \rbrace$ for all $\delta >0$.

Assume that $NZ_{\delta, > }(r)(A)$ is finite for all $r$ and $\delta $. Then $A$ is $\sigma$-finite for $\mu$.
\end{lemma}

\vspace{2 mm}

\underline{\textbf{Demonstration:}} We will also use another notation: $NZ_{\delta , \geq }(r)(A) = \lbrace z \in L / \mu ((A-r)\cap (C+z))\geq \delta \rbrace$. Fist, we prove that $NZ_{\delta, \geq }(.)(A): C \rightarrow P(L)$ is upper semi-continuous in weak topology for any strong Borel set $A$ and any $\delta > 0$. Let $(u_n)$ be a sequence of $C$ weakly converging to a $u_{\infty}$. Let us show that if $z\in NZ_{\delta, \geq }(u_n)(A)$ for $n$ tall enough, then $z\in NZ_{\delta, \geq }(u_{\infty})(A)$. That is, if $\mu ((A-u_n)\cap (C+z))\geq \delta$ for $n$ tall enough, $\mu ((A-u_{\infty})\cap (C+z))\geq \delta$. Choose a $\varepsilon \in (0;\delta )$, and let $O$ be a weak open of $C+z$ such that $(A-u_{\infty})\cap (C+z)) \subset O$ and $\mu (O) < \mu ((A-u_{\infty})\cap (C+z))) + \varepsilon$. For all $n$, let $F_n$ be a weak closed included in $(A-u_n)\cap (C+z)$ with $\mu (F_n) \geq \delta - \varepsilon$. Since $C+z$ is a weak compact, the set of weak closed of $C+z$ is compact for the weak Hausdorff distance. Hence, one can choose a subsequence of $(F_n)$ converging to a $F$, and it is clear that $F \subset O$. Hence $\mu (O) \geq \mu (F) \geq \limsup_{n \rightarrow +\infty} \mu (F_n) \geq\delta - \varepsilon$. Indeed $(F_n)$ approaches $F$ in weak Hausdorff distance, and $\mu$ is regular for weak topology. So: $\mu((A-u_{\infty})\cap (C+z))) \geq\delta - 2\varepsilon$ for all $\varepsilon > 0$. Therefore $\mu ((A-u_{\infty})\cap (C+z))) \geq\delta$.

Hence $NZ_{\delta, \geq }(.)(A): C \rightarrow P(L)$ is upper semi-continuous in weak topology.

\vspace{2 mm} 

Now, note that $z \in NZ_{\delta, > }(r)(A) \Leftrightarrow \mu ((A-r)\cap (C+z))> \delta$, hence $z \in NZ_{\delta, > }(r)(A) \Leftrightarrow \mu (((E\setminus A)-r)\cap (C+z))< 1-\delta$, so $z \in NZ_{\delta, > }(r)(A) \Leftrightarrow z \not\in NZ_{1-\delta, \geq }(r)(E\setminus A)$. Consequently, $L\setminus NZ_{\delta, > }(r)(A) = NZ_{1-\delta, \geq }(r)(E\setminus A)$ for all $\delta \in (0,1)$, $r \in C$, and all strong Borel set $A$.

\vspace{2 mm}

We will now use a compactness argument to show that $A$ is $\sigma$-finite. Let $F$ be any finite subset of $L$, and denote $C_F = \lbrace r \in C/ NZ_{\delta, > }(r)(A) \subset F \rbrace $, assuming that $NZ_{\delta, > }(r)(A)$ is finite for all $r$ and $\delta $. It is obvious that the family $(C_F)_{F \mbox{ finite}}$ covers $C$. Moreover, $C$ is compact in weak topology. Let us show that $C_F$ is a weak open for any finite $F$. 

Indeed, $C\setminus C_F = \lbrace r \in C / NZ_{\delta, > }(r)(A) \not\subset F\rbrace = \lbrace r \in C / NZ_{1-\delta, \geq }(r)(E\setminus A) \cap F \neq \emptyset\rbrace$. Now, we will show that $\lbrace r \in C / NZ_{1-\delta, \geq }(r)(E\setminus A) \cap F \neq \emptyset\rbrace$ is a weak closed for any $\delta \in (0,1)$. Let $(r_n)$ be a sequence of that set converging to, say, $r_{\infty }$. For any $n$, we can choose a $z_n \in NZ_{1-\delta, \geq }(r_n)(E\setminus A) \cap F$, and since $F$ is finite, we can find a strictly increasing sequence $\varphi $ of natural numbers such that $(z_{\varphi (n)})$ is constant, equalling, say, $z_{\infty } \in F$. Then $(r_{\varphi (n)})$ approaches $r_{\infty }$, and $z_{\infty } \in NZ_{1-\delta, \geq }(r_{\varphi (n)})(E\setminus A) \cap F$ for all $n$. Hence, since $NZ_{\delta, \geq }(.)(B): C \rightarrow P(L)$ is upper semi-continuous in weak topology for any strong Borel set $B$, $z_{\infty } \in NZ_{1-\delta, \geq }(r_{\infty })(E\setminus A) \cap F$, and $C\setminus C_F $ is a weak closed. Consequently, $C_F$ is a weak open for any finite set $F \subset L$.

\vspace{2 mm}

So, $C$ is a weak compact covered by the family of opens $(C_F)_{F \mbox{ finite}}$. Hence, it can be covered by only a finite number of them, say $(C_{F_i})$, and so $C \subset C_{\bigcup_i F_i}$. Then, there is a finite set $\Phi_{\delta }(A) = \bigcup_i F_i \subset L$ such that $C = C_{\Phi_{\delta }(A)}$. Consequently, for any $r \in C$, $ NZ_{\delta, > }(r)(A) \subset \Phi_{\delta }(A) $.

\vspace{2 mm}

Now, since $C$ is a weak compact, there is a denumerable weakly dense subset $P$ in $C$. And $M_{\delta}(A)= \bigcup_{r \in P} \bigcup_{z \in \Phi_{\delta }(A)} z+r+C$ is a denumerable union of strong $1/2$-balls. Now, for any $1/2$-ball $B$, $m(A\cap B \setminus M_{\delta }(A))\leq \delta$. Let us prove this.

Let $B$ be any $1/2$-ball. We can write $B = z+r+C$ with $z \in L$ and $r \in C$. Then $r$ is a weak limit of a sequence of $P$, say $(r_n)$. For any $n$, $m(A\cap (z+r_n+C) \setminus M_{\delta }(A))\leq \delta$. Indeed, if $z \in NZ_{\delta, > }(r_n)(A)$, $m(A\cap(z+r_n+C)\setminus M_{\delta}(A))$ vanishes, and if $z \not\in NZ_{\delta, > }(r_n)(A)$, $m(A\cap (z+r_n+C) \setminus M_{\delta }(A))\leq m(A\cap (z+r_n+C)) \leq \delta$. Hence, for any $n$, $z \not\in NZ_{\delta, > }(r_n)(A\setminus M_{\delta }(A))$, that is, $z \in NZ_{1-\delta, \geq }(r_n)(E\setminus (A\setminus M_{\delta }(A)))$. Since $(r_n)$ weakly converges to $r$, $z \in NZ_{1-\delta, \geq }(r)(E\setminus (A\setminus M_{\delta }(A)))$, and $z \not\in NZ_{\delta, > }(r)(A\setminus M_{\delta }(A))$. Consequently, $m(A\cap (z+r+C) \setminus M_{\delta }(A))\leq \delta$.

\vspace{2 mm}

Then, denote $M = \bigcup_{n \in \mathbf{N}} M_{1/n }(A)$. It is a denumerable union of denumerable unions of finite sets of $1/2$-balls, so it is a countable union of $1/2$-balls. And for any strong $1/2$-ball $B$, $m(A\cap B \setminus M)\leq m(A\cap B \setminus M_{\delta }(A)) < \delta$ for all $\delta > 0$. Then, $m(A\cap B \setminus M) = 0$, and $A \setminus M$ is null for $1/2$-balls, hence null. Consequently, $A$ is $\sigma $-finite.\begin{flushright}
$\blacksquare$
\end{flushright}

As a consequence, we see that if a Borel set $A$ has finite measure, it is $\sigma$-finite.

\subsection{A technique for a easy calculus of integrals}

We will now describe a easy way to calculate integrals for $\mu$ in $E = l^{\infty}$.

First of all, we need to know some properties of integrable functions. Let $f : l^{\infty} \rightarrow \mathbf{R}$ be any measurable function, $(x_i)$ the canonical basis of $l^{\infty}$, $L$ the associated lattice and $C$ its unit cell. It is well known that $f$ can be written $f = f^{+}-f^{-}$ with $f^{+}$ and $f^{-}$ measurable non-negative functions such that $\vert f\vert = f^{+}+f^{-}$. $f^{+}$ and $f^{-}$ are called the non-negative and the non-positive parts of $f$. $f$ is integrable if and only if $f^{+}$ and $f^{-}$ are integrable.

\subsubsection{Notion of support of a measurable function}

We begin with recalling an elementary fact dealing with the support of a function and its integral. We use here the notion of measurable support which is suitable in measure theory. Usually, the support of a function is the closure of the measurable support defined here.

\vspace{2 mm}

\begin{lemma} 
Let $f$ be any measurable non-negative function, and $A$ any measurable set for any measure $m$. $\int_A fd\mu = 0$ if and only if $m(\lbrace x \in A / f(x) \neq 0\rbrace)=0$. 
\end{lemma}

\vspace{2 mm}

\underline{\textbf{Demonstration:}} We have: $\lbrace x\in A / f(x) \neq 0\rbrace = \bigcup_{n \in \mathbf{N}^*} \lbrace x / f(x) \geq 1/n\rbrace$. Hence, if $m(\lbrace x\in A / f(x) \neq 0\rbrace)\neq 0$, there is a $n \in \mathbf{N}^*$ such that $m(\lbrace x / f(x) \geq 1/n\rbrace)\neq 0$, and $\int_A fd\mu \geq \int_{\lbrace x / f(x) \geq 1/n\rbrace} fd\mu \geq \frac{1}{n}m(\lbrace x / f(x) \geq 1/n\rbrace) >0$, which is absurd. The converse is obvious.
\begin{flushright}
$\blacksquare$
\end{flushright}

We use the following definition:

\vspace{2 mm}

\begin{definition} 
Let $f$ be any measurable function on $E=l^{\infty}$. The set $\lbrace x \in A / f(x) \neq 0\rbrace$ is called the (measurable) support of $f$ and denoted $Supp(f)$.
\end{definition}

\vspace{2 mm}

The support of a measurable function is obviously a measurable set.

\subsubsection{Integrability of measurable functions}

We will now give a necessary condition for a non-negative function to be integrable using that notion of support.

\begin{proposition} 
Let $f$ be any non-negative measurable function on $E=l^{\infty}$. Assume that $f$ is integrable. Then $Supp(f)$ is $\sigma$-finite.
\end{proposition}

\vspace{2 mm}

\underline{\textbf{Demonstration:}} Assume that $\int_E fd\mu < + \infty$. Choose a $\varepsilon >0$, and assume that $Supp(\max (f-\varepsilon,0))$ is not $\sigma$-finite. Then $\int_{Supp(\max (f-\varepsilon,0))} fd\mu \leq \int_E fd\mu < + \infty$, and $\int_{Supp(\max (f-\varepsilon,0))} fd\mu \geq \varepsilon \mu (Supp(\max (f-\varepsilon,0))) = + \infty $. It is absurd. Hence $Supp(\max (f-\varepsilon,0))$ is $\sigma$-finite for all $\varepsilon > 0$, and $Supp(f) = \bigcup_{n \in \mathbf{N}} Supp(\max (f-1/n,0))$ is a denumerable union of countable unions of $1/2$-balls with null sets, hence is $\sigma$-finite. \begin{flushright}
$\blacksquare$
\end{flushright}

Then we have :

\begin{theorem}[Characterisation of integrable functions] 
Let $f$ be any measurable function on $E$. $f$ is integrable if and only if both its non-negative and non-positive parts are.

Let $f^+$ be a non-negative function. It is integrable if and only if its support is $\sigma$-finite for an almost everywhere countable covering by $1/2$-balls $(B_i)$, and if the family $(\int_{B_i \setminus \bigcup_{j<i}B_j} fd\mu)_i$ is a convergent series.

Let $f$ be a measurable function on $E$. It is integrable if and only if its support is $\sigma$-finite for an almost everywhere countable covering by $1/2$-balls $(B_i)$, if for all $i$, $f$ is integrable on $B_i \setminus \bigcup_{j<i}B_j$ and if the family $(\int_{B_i \setminus \bigcup_{j<i}B_j} \vert f\vert d\mu)_i$ is a convergent series.
\end{theorem}

\vspace{2 mm}

\underline{\textbf{Demonstration:}} The first point is well known. The second is quite obvious with what we showed. For the third one, note that $Supp(f) = Supp(f^+)\cup Supp(f^-)$, so $Supp(f)$ is $\sigma$-finite if and only if both $Supp(f^+)$ and $Supp(f^-)$ are. Use then the characterisation given by the two first points and the well known results about integrability and convergence of series.
\begin{flushright}
$\blacksquare$
\end{flushright}

\subsubsection{Integrals of integrable functions}

First of all, we have the obvious :

\begin{proposition} 
Let $f$ be any integrable function on $E=l^{\infty}$. Then : 

$$\int_E f d\mu = \sum_{i} \int_{B_i \setminus \bigcup_{j<i}B_j} f d\mu ,$$

where $(B_i)$ is an almost everywhere countable covering of $Supp(f)$ by $1/2$-balls.
\end{proposition}

\vspace{2 mm}

We need now to find an easy way to calculate integrals of integrable functions on $1/2$-balls.

Let $f$ be such an integrable function defined on a parallelepiped $P$ containing a strong $1/2$-ball $C$. Denote $(c_i)_i$ the center of $C$, $d_i = c_i-1/2$ for all $i$, and let $(a_i)_i$ be any element of $P$. For example, one can choose $a_i=0$ for all $i$ if $0 \in P$. $f$ can be seen as a function of a sequence of numbers of $[0;1]$ (in the canonical basis, recall that this basis is in reality a weak Schauder basis): 

$$f((x_i)_i)=f(d_0 + x_0;d_1+x_1;d_2+x_2;...;d_n+x_n;....). $$

Denote $f_n : [0;1]^{\mathbf{N}} \rightarrow \mathbf{R}$ the map such that :

$$f_n ((x_i)_i) = f(d_0+x_0;d_1+x_1;d_2+x_2;...;d_{n-1}+x_{n-1};d_n+x_n;a_{n+1};a_{n+2};...).$$

This definition makes sense since $P$ is a parallelepiped. It is clear that for any sequence of numbers $(x_i)$ of $[0;1]^{\mathbf{N}}$, the sequence $(\sum_{i<n} (d_i +x_i ) e_i + \sum_{i\geq n} a_i e_i)_n$ approaches $\sum_{i} x_i e_i$ in weak topology ($(e_i)$ denotes the canonical basis). Moreover, since $f$ is integrable on $C$ which is weakly metric compact, with Luzin Theorem, there is a null set $N$ such that $f$ is weakly continuous out of $N$. Hence, $(f_n)$ simply converges to $f$ out of $N$. Suppose now that $f$ is bounded by, say, $M$, on $P$, except on a null set. Then, with Lebesgue convergence theorem :

$$\lim_{n\rightarrow +\infty} \int_C f_n d\mu =\int_C fd\mu .$$

More generally, without any assumption about $f$, for all $M$,

$$\lim_{n\rightarrow +\infty} \int_C f_n (x)\textbf{1}_{\lbrace \vert f_n (x)\vert \leq M\rbrace} d\mu = \int_C f (x)\textbf{1}_{\lbrace \vert f (x)\vert \leq M\rbrace} d\mu ,$$

and :

$$\lim_{M\rightarrow +\infty}\lim_{n\rightarrow +\infty} \int_C f_n (x)\textbf{1}_{\lbrace \vert f_n (x)\vert \leq M\rbrace} d\mu = \int_C fd\mu.$$

Indeed, $fd\mu$ is a finite non-negative measure, and $\bigcap_{n \in \mathbf{N}} \lbrace \vert f(x)\vert > n\rbrace = \emptyset$ since $f : P \rightarrow \mathbf{R}.$

Now, in order to easily calculate integrals like $\int_C f_n d\mu $, we will use a simple form of Fubini Theorem for $\mu$ on $C$. $\mu$ equals on $C$ the measure given by the Geometry and Probability Theorem for any limit notion (by uniqueness) because of the patching measures Theorem. And the Geometry and Probability Theorem states that once a ultra-filter is given, for any compact spaces $K$ and $K'$, $\mu_K \otimes \mu_{K'} = \mu_{K\times K'}$. Here, for any $n \in \mathbf{N}$, $[0;1]^{\mathbf{N}}$ can be written as the product of $[0;1]^{n+1}$ for the n first coordinates with $[0;1]^{\lbrace i> n+1, i \in \mathbf{N}\rbrace}$ for the last ones (they are in denumerable cardinal), and both these spaces are compact in weak topology. Hence, for $n\geq 1$ : 

$$\int_{C} f_n d\mu = \int_{[0;1]^{n+1}} (\int_{[0;1]^{\lbrace i> n+1, i \in \mathbf{N}\rbrace }}f_n d\mu_{[0;1]^{\lbrace i> n+1, i \in \mathbf{N}\rbrace }})d\mu_{[0;1]^{n+1}} $$
$$= \int_{[0;1]^{n+1}} f_n\vert_{[0;1]^{n+1}\times \lbrace (a_i)_{i \geq n+1}\rbrace} d\lambda.$$

So : 

$$ \lim_{n \rightarrow +\infty}\int_{[0;1]^{n+1}} f_n\vert_{[0;1]^{n+1}\times \lbrace (a_i)_{i \geq n+1}\rbrace} d\lambda. = \int_{C} f d\mu ,$$

for any bounded function. 
We showed :

\begin{theorem} 
Let $f$ be any integrable function on a parallelepiped containing a given element $(a_i)_i$ and a $1/2$-ball, say $(d_i)_i + [0,1]^{\mathbf{N}}$. Denote $f_n ((x_i)_{i\leq n}) = f(d_0 + x_0;d_1 + x_1;d_2 + x_2;...;d_{n-1} + x_{n-1};d_n + x_n;a_{n+1};a_{n+2};...)$. 

- Assume that $f$ is bounded on $P$. Then the limit $\lim_{n \rightarrow +\infty}\int_{[0,1]^{n+1}} f_n d\lambda$ exists, and we have :

$$ \lim_{n \rightarrow +\infty}\int_{[0,1]^{n+1}} f_n d\lambda = \int_{C} f d\mu ,$$

where $d\lambda$ is the Lebesgue measure of finite dimensional spaces.

- Without any more assumption about $f$,

$$\int_{C} f d\mu = $$
$$\lim_{M\rightarrow +\infty} \lim_{n \rightarrow +\infty} \int_{[0,1]^{n+1}} f_n \textbf{1}_{\lbrace \vert f_n\vert \leq M\rbrace}d\lambda.$$

\end{theorem}

\vspace{2 mm}

Moreover, using Luzin and Lebesgue Theorems, we can show :

\begin{proposition} 
Let $f$ be any measurable function on a parallelepiped containing a given element $(a_i)_i$ and a $1/2$-ball, say $(d_i)_i + [0,1]^{\mathbf{N}}$. We use the same notations as above. Then $f$ is integrable if and only if all $M$, $\lim_{n \rightarrow +\infty}\int_{[0,1]^{n+1}} \vert f_n\vert\textbf{1}_{\lbrace \vert f_n\vert\leq M\rbrace} d\lambda $ exist, and: 
$$\lim_{M\rightarrow +\infty} \lim_{n \rightarrow +\infty}\int_{[0,1]^{n+1}} \vert f_n \vert\textbf{1}_{\lbrace \vert f_n \vert\leq M\rbrace} d\lambda < +\infty .$$
\end{proposition}

Beware to the fact that one can have $\lim_{n\rightarrow +\infty} \int_{[0;1]^{n+1}} f_n d\lambda \neq \int_C fd\mu$ when $f$ is integrable but unbounded. For instance, look at the following function :

\begin{tabbing}
$f :$ \= $[0;1]^{\mathbf{N}}$\= $\rightarrow$ \= $\mathbf{R} $\\
\> $(x_i)_i$\> $\mapsto$ \> $\frac{3}{2}\textbf{1}_{\lbrace x_0 \in [0;1/3]\cup [2/3;1]\rbrace}\times\prod_{i \geq 1}\textbf{1}_{\lbrace x_i \in [0;1/3]\rbrace}+$
\end{tabbing}

$$\sum_{n=1}^{\infty} \frac{3^{n+1}-3^n}{2^n}\prod_{i<n} \textbf{1}_{\lbrace x_i \in [0;1/3]\cup [2/3;1]\rbrace}\times \textbf{1}_{\lbrace x_n \in [2/3;1]\rbrace}\times \prod_{i>n} \textbf{1}_{\lbrace x_i \in [0;1/3]\rbrace}.$$

$f$ is measurable for the strong $\sigma$-algebra, and $Supp(f) = ([0;1/3]\cup [2/3;1])^{\mathbf{N}}$. Moreover, 

$$\mu (Supp(f)) = \int_{Supp(f)} d\mu = \lim_{n \rightarrow +\infty } \int_{([0;1/3]\cup [2/3;1])^{n+1}} d\lambda$$
$$=\lim_{n \rightarrow +\infty } (\frac{2}{3})^{n+1}=0,$$
since $x \mapsto 1$ is bounded. Hence, $\int_{[0;1]^{\mathbf{N}}} fd\mu = 0$. 

But:
$$\int_{[0;1]^{N+1}} f(x_0;...;x_N;0;...;0;...)d\lambda =$$
$$2\times \frac{3}{2}(\frac{1}{3})^{N+1}+ \sum_{n=1}^{N}\frac{3^{n+1}-3^n}{2^n}(\frac{2}{3})^n (\frac{1}{3})^{N-n+1} =$$
$$\frac{1}{3^{N+1}}(\sum_{n=1}^{N} (3^{n+1}-3^n)+3)=1$$ for all $N$. 

\vspace{10 mm}

We can now give a method to calculate integrals for $\mu$.

\vspace{2 mm}

Let $f$ be any measurable function on $E$, and denote $\bar{f}_n((x_i)_{i\leq n}) = f(x_0;x_1;...;x_{n-1};x_n;0;0;...)$.

- Calculate $Supp(f)$.

- If $Supp(f)$ is not $\sigma$-finite, $f$ is not integrable. If it is, choose an almost everywhere countable covering of $Supp(f)$ by $1/2$-balls $(B_i)$.

- Calculate $ \int_{B_i \setminus \bigcup_{j<i}B_j} \vert f\vert d\mu$ for all $i$ with the formula :

$$  \int_{B_i \setminus \bigcup_{j<i}B_j} \vert f\vert d\mu =$$

$$\lim_{M \rightarrow +\infty} \lim_{n \rightarrow +\infty}\int_{(B_i \setminus \bigcup_{j<i}B_j)\cap \mathbf{R}^{n+1}} \vert \bar{f}_n \vert \textbf{1}_{\lbrace \vert \bar{f}_n \vert\leq M\rbrace}d\lambda (x_0;...;x_n).$$

- If some of them are infinite, $f$ is not integrable. Otherwise,

- Calculate :

$$\sum_{i} \int_{B_i \setminus \bigcup_{j<i}B_j} \vert f\vert d\mu .$$

If it is infinite, $f$ is not integrable, and if not, $f$ is integrable.

- Calculate $ \int_{B_i \setminus \bigcup_{j<i}B_j} f d\mu$ for all $i$ with the formula :

$$  \int_{B_i \setminus \bigcup_{j<i}B_j} f d\mu =$$
$$\lim_{M\rightarrow + \infty}\lim_{n \rightarrow +\infty}\int_{(B_i \setminus \bigcup_{j<i}B_j)\cap \mathbf{R}^{n+1}} \bar{f}_n \textbf{1}_{\lbrace \vert \bar{f}_n \vert \leq M\rbrace} d\lambda (x_0;...;x_n).$$

- Finish with :

$$\int_E fd\mu = \sum_{i} \int_{B_i \setminus \bigcup_{j<i}B_j}  f d\mu .$$

\vspace{2 mm}

Beware to the fact that $f$ can have a $\sigma$-finite support, with $f$ integrable on all $B_i \setminus \bigcup_{j<i}B_j$ (where $(B_i)_i$ is a almost everywhere covering of $Supp(f)$), and such that the sum $\sum_{i} \int_{B_i \setminus \bigcup_{j<i}B_j}  f d\mu$ converges, but without integrability of $f$. In this case, $\sum_{i} \int_{B_i \setminus \bigcup_{j<i}B_j} \vert f\vert d\mu = + \infty $.

\subsection{Invariance by translations}

With the previous calculus method and invariance of $\mu$ by translations, we can easily show :

\begin{theorem} 
Let $f$ be any function on $E$, and $t \in E$. $f$ is integrable if and only if $x \mapsto f(x+t)$ is, and in this case, 

$$\int_E fd\mu = \int_E f(x+t)d\mu (x).$$

\end{theorem}

\subsection{Fubini Theorems}

By the same method, we can show :

\begin{theorem} 
Let $ E = l^{\infty}$, $(e_i)$ be its canonical basis. Let $V$ be a finite dimensional vector space of $E$ admitting a basis made of elements of $(e_i)$, and $W$ be the strong closure of the space generated by the other elements of $(e_i)$. We have : $E = V\oplus W$. Then, for any integrable function $f$,

$$\int_E fd\mu = \int_V (\int_W fd\mu_W)d\lambda_V,$$

where $d\mu_W$ is the previously built measure associated to $W$, a space isomorphic to $E$, and $d\lambda_V$ the usual Lebesgue measure of $V$. Moreover, any function $g : E \rightarrow \mathbf{R}$ is integrable if and only if $g\vert_{W+\lbrace x\rbrace}$ is integrable for $\lambda$-almost all $x \in V$, and if $x\mapsto \int_W g\vert_{W+\lbrace x\rbrace} d\mu_W$ is integrable on $V$.
\end{theorem}

The proof is straightforward with the usual Fubini Theorem and the behaviour of $\mu$ with products.

\vspace{2 mm}

We are sure that a generalisation of the previous result can be made weakening the assumption that $V$ and $W$ are built with the canonical basis. It demands a geometric theory of $l^{\infty }$ using the fact that the canonical basis is a weak Schauder basis. It has to be linked with a notion of Jacobian determinant, and with a notion of linear continuous maps preserving the measure, as in the finite dimensional case.

\vspace{2 mm}

Nevertheless, something quite new appears in our frame of work.

\begin{theorem} 
Let $ E = l^{\infty}$, and $(e_i)$ be its canonical basis. Let $V'$ be an infinite dimensional vector space of $E$ admitting a basis made of elements of $(e_i)$, but such that infinitely many $e_i$ do not belong to $V'$, $V$ its strong closure, and $W$ be the strong closure of the space generated by the other elements of $(e_i)$. We have : $E = V\oplus W$. Then, for any integrable function $f$,

$$\int_E fd\mu = \int_V (\int_W fd\mu_W)d\mu_V,$$

where $d\mu_W$ and $d\mu_V$ are the previously built measure associated to $W$ and $V$. Moreover, any function $g : E\rightarrow \mathbf{R}$ is integrable if and only if $g\vert_{W+\lbrace x\rbrace}$ is integrable for $\mu_V$-almost all $x \in V$, and if $x\mapsto \int_W g\vert_{W+\lbrace x\rbrace} d\mu_W$ is integrable on $V$.
\end{theorem}

The proof follows the same way as before.

\vspace{2 mm}

These theorems suggest that a theorem of the following kind holds, with good hypothesis.

\vspace{2 mm}

Assume that $E = \bigoplus_i V_i$, a denumerable direct sum of closed vector subspaces of $E$, such that for any $i$, $V_i$ is the strong closure of its subspace generated by the elements of the canonical basis contained in. This assumption contains the geometry, and can surely be weakened. Then, for any bounded integrable function $f$, $ \lim_{n\rightarrow + \infty}\int_ {V_0} (\int_{V_1}...(\int_{V_n} fd\mu_{V_n})...d\mu_{V_1})d\mu_{V_0}$ exists, and

$$\int_E fd\mu = \lim_{n\rightarrow + \infty}\int_ {V_0} (\int_{V_1}...(\int_{V_n} fd\mu_{V_n})...d\mu_{V_1})d\mu_{V_0},$$

where $d\mu_{V_i}$ can be a measure of the previous type or a Lebesgue measure on a finite dimensional space for all $i$. Maybe some other analytic assumption has to be made about $f$ (probably something about the support of $f$ and null sets).

\section*{Conclusion}

\vspace{15 mm}

{\Large Greeting} I thank M. Paulin (Universit\'{e} Paris-Sud), M. Penot
(Universit\'{e} des pays de l'Adour) and M. Ciligot-Travain (Universit\'{e}
d'Avignon) for the interest they manifested for my work. I express my
gratitude to Mme Kosmann-Schwarzbach, M. Aubin (Professeurs honoraires),
M. Zvonkine (Universit\'{e} de Bordeaux I) and specially to M. Lewandowski (Université de Rennes) and M. Léonard (Université Paris-Ouest) for their encouragements and their help. Finally, I
warmly thank M. Drouin and M. Paya, Teacher like me at Lyc\'{e}e de Borda in Dax, for
their work of re-lecture and check of demonstrations, and for our fruitful discussions.

\end{document}